# Congestion Management by Applying Co-operative FACTS and DR program to Maximize Renewables

Farhad Samadi Gazijahani[1*], Rasoul Esmaeilzadeh[2]
[1,2] Azarbaijan Regional Electric Company, Tabriz, Iran
Email address: [1] farhad.samadi.g@gmail.com, [2] rasoul_zadeh@yahoo.com

**Abstract:** This research proposes an incremental welfare consensus method based on flexible alternating current transmission systems (FACTS) and demand response (DR) programs to control transmission network congestion in order to increase the penetration of wind power. The locational marginal prices are used as an input by the suggested model to control the FACTS device and DR resources. In order to do this, a cutting-edge two-stage market clearing system is created. In the first stage, participants bid on the market with the intention of maximizing their profits, and the ISO clears the market with the goal of promoting societal welfare. The second step involves the execution of a generation re-dispatch issue in which incentive-based DR and FACTS device controllers are optimally coordinated to reduce the rescheduling expenses for generating firms. Here, a static synchronous compensator and a series capacitor operated by a thyristor are used as two different forms of FACTS devices. A case study on the modified IEEE one-area 24-bus RTS system is then completed. The simulation results show that the suggested interactive DR and FACTS model not only reduces system congestion, but also makes the system more flexible so that it can capture as much wind energy as feasible.

**Keywords:** FACTS; DR program; Renewable harvesting; Congestion management.

## 1. Introduction

Following power system reorganization, there has been a sharp rise in the inclination for the many independent participants in the energy market to maximize their profits, which has led to intense use of the transmission networks [1]. Due to this problem, the power systems are often used at or close to their maximum capacity. The advent of any unforeseen contingency in the transmission network or an abrupt rise in load demand will result in overloading and congestion in some or even all network lines under these conditions [2].

Because of network limitations, congestion not only raises the price of providing electricity but also reduces the system's ability to generate enough income [3].

On the other hand, during the last ten years, the power systems have gradually been infiltrated by renewable resources, mainly wind turbines (WT), which has led to the creation of new issues in the safe and dependable operation of power systems. Additionally, the intermittent nature and intrinsic fluctuation of wind energy have made this situation more complicated and added substantial challenges [2]. In light of this, the independent system operator (ISO) should take the necessary steps to reduce and control the congestion of wind integrated transmission networks, which if left unchecked may raise locational marginal pricing (LMP) and jeopardize the safety and stability of the electricity systems.

The term "congestion management problem" (CMP) refers to a series of procedures for effectively reducing transmission network congestion while remaining within the bounds of the system [4]. These procedures may be divided into two primary categories, consisting of cost-free and non-cost-free techniques. Flexible alternating current transmission systems (FACTS), a cutting-edge tool for compensating reactive power, may now increase line capacity by reducing reactance on crowded lines [5]. Additionally, the importance of demand response (DR) programs in reducing such congestions is growing as a result of the introduction of advanced measuring infrastructure (AMI) in the power sector [6]. However, there hasn't been much work put towards collaborating the DR program and FACTS on LMP smoothing and congestion relief. A appropriate solution that fulfills the system's operational needs should also be used to account for the unpredictability associated with wind generating. By strengthening the cooperation between the ISO and FACTS in executing DR programs, this research presents a probabilistic incremental welfare consensus model for CMP of wind integrated transmission networks to fill this gap in the literature.

In the technical literature, the CMP has recently been the subject of considerable investigation in a deregulated context. Reactive power support, generation rescheduling strategies, operation of transformer taps, and load shedding are some of the common CMP procedures [4]. The active power of GENCOs was optimally rescheduled in [7] using a genetically based optimization algorithm in order to reduce transmission line congestion and increase the social welfare of the system by reducing the difference in LMP across buses. In [8], a decentralized approach of decision-making for multi-utility transmission networks was put out. In order to discover the optimum overall solution with the goal of improving the system's social welfare, an interior point based optimization technique was provided in this work. Additionally, the decentralized approach's findings were compared to and validated using the centralized method. The best way to use renewable energy sources (RESs) to reduce traffic in a deregulated environment was researched in [9]. Based on an apparent congestion index, the suggested solution reduces the cost of re-dispatching both renewable and traditional energy sources. With the aid of power electronic interfaces, the static synchronous compensator (STATCOM) also supports reactive power while rescheduling the reactive power of RES. In order to prevent network line congestion, the best participating generators selection method was described in [10]. Additionally, the particle swarm optimization (PSO) technique was used to schedule GENCOs more effectively, and the most efficient generator for overcoming congestion was chosen based on how sensitive each generator was to the power flows coming from congested lines.

Building intelligent electricity systems provides adequate incentive for engaged customers to contribute to boosting system efficiency. Actually, the DR program, which is a key CMP tool, makes it possible to encourage clients to take part in CMP tasks that attempt to smooth LMPs. Choosing the best location and timing for DR deployment was suggested in [11], which reduced network congestion by reducing superfluous demands. As congestion criteria, power transfer distribution factor and available transfer capability (ATC) are used to identify the best locations for DR resources. In [12], the price area approach was established to reduce transmission line congestion. For two and three area systems, this notion and associated control approach are examined. By raising prices in deficit areas to address the congestion issue, more will be produced while less will be consumed. Additionally, a novel CMP strategy using a multi-objective optimization algorithm linked to mechanisms for load shedding and generation re-dispatching was put out in [13]. The objective function additionally takes demand bidding into account and includes load served error, generation and load shedding cost reduction, load served and social welfare maximization. The use of FACTS devices to improve the technical properties of renewable-integrated power systems has received a lot of attention recently. For instance, [14] described an integrated CMP technique based on FACTS controllers and generator rescheduling. According to the findings, using FACTS compensators lowers the cost of re-dispatching generators. A bi-level optimization approach was given by the authors in [15] to control network overloading by selecting a certain DR contract in a day-ahead market. The advised bi-level approach resulted in accurate ATC. To increase the resilience and participation levels of a wholesale market in a real network, a novel model was also introduced in [16]. The experimental strategy being suggested is to establish a trend that will include more clients in DR programs. The authors of [17] demonstrated two different kinds of FACTS devices, an optimum parameter setting for power system variables, and a novel fuzzy logic membership function for identifying weak nodes. The major goal of this work was to determine the best position for FACTS controllers in order to minimize power losses and ease transmission line congestion. Using dispatchable generating and responsive loads, the authors of [18] demonstrated a distributed solution to DC optimum power flow (OPF) for coordinating all participants to optimize social welfare while reducing transmission line congestion. Additionally, storage and electric car proposals were made in [19] to reduce the load on active distribution networks while taking into account heterogeneous distributed generation (DG) resources.

The balancing of reactive power across subsystems is one of the most important fundamental concerns in the CMP of power systems. In a reorganized hybrid power market, an optimum CMP with reactive power assistance was created in [20]. Rearranging both the active and reactive power generating schedules was also employed to control line congestion. This paper's major goal was to reduce the cost of rescheduling using a PSO-based approach. To simultaneously maximize congestion cost, voltage security, and transient security, a multi-objective CMP was recommended in [21]. To choose the best effective outcome from the Pareto frontier obtained from the multi-objective optimization, the authors used a fuzzy satisfactory technique. The best size and position of several FACTS devices in a large-scale transmission network were determined using a graphical user interface-based model that relied on a genetic algorithm in [22] to enhance power system loadability. The best positioning of the thyristor controlled series capacitor (TCSC) for reducing congestion and boosting ATC was explored in [23]–[25]. For instance, a sensitivity factor-based model that allocated TCSC in the distribution network with the goal of reducing both congestion and investment costs was successfully implemented in [23]. To increase the voltage and transient stability margins while lowering the system's overall operating cost, multi-objective congestion management was utilized in [24]. The authors of [25] provided a technique for placing TCSCs that uses the congestion rent contribution method and LMP calculation.

The current study aims to make the following contributions in order to address the flaws and shortcomings of earlier works:

- Creating an integrated CMP architecture using LMP to effectively manage the DR and FACTS.
- To enhance social welfare and reduce network congestion at the lowest feasible cost, an unique two-stage market clearing method is devised.
- Making use of LMP signaling, which ISO transmits to FACTS devices, to cooperatively modify the set points of FACTS devices in accordance with DR aggregators for congestion reduction.
- Deciding how best to use DR resources to help FACTS get rid of congestion.

The remainder of the essay is structured as follows. The description of the DR program and the modeling of FACTS devices are described in parts 2 and 3, respectively. Section 4 describes a two-stage market clearing process and a formulation for congestion management. Section 5 provides the simulation findings from a real-world case study. Section 6 closes the essay and identifies some positive aspects.

## 2. FACTS devices modeling

There is no way to complete all contracts between consumers and producers in the restructured environment owing to the emergence of congestion in the transmission lines. This problem will result in a rise in LMP at various nodes, which lowers the system's social welfare. In order to address this problem, FACTS devices, particularly series FACTS devices, may be effectively used to reduce power flows in densely loaded lines, resulting in LMP smoothing, while also lowering power losses and procurement costs [14]. The static model of these controllers is used in the current research to offer reactive compensation utilizing FACTS devices for improving the capability of transmission lines. The FACTS device is taken into account by the static power injection model as a component that injects active and reactive power at terminal buses with the device positioned between them.

## 2.1 TCSC

The TCSC can be modeled as a capacitive reactance $(-jx_{TCSC})$ during the steady-state condition [26]. The change in the line flow due to the mentioned series compensation can be reflected as power injections at terminal buses $i$ and $j$. Consequently, the active and reactive powers injections of TCSC at buses $i$ and $j$ are given as follows [26]:

$$P_i^C = |V_i|^2 \Delta G_{ij} - |V_i||V_j|[\Delta G_{ij} \cos \delta_{ij} + \Delta B_{ij} \sin \delta_{ij}] \tag{1}$$

$$Q_i^C = -|V_i|^2 \Delta B_{ij} - |V_i||V_j|[\Delta G_{ij} \sin \delta_{ij} - \Delta B_{ij} \cos \delta_{ij}] \tag{2}$$

$$P_j^C = |V_j|^2 \Delta G_{ij} - |V_i||V_j|[\Delta G_{ij} \cos \delta_{ij} - \Delta B_{ij} \sin \delta_{ij}] \tag{3}$$

$$Q_j^C = -|V_j|^2 \Delta B_{ij} + |V_i||V_j|[\Delta G_{ij} \sin \delta_{ij} + \Delta B_{ij} \cos \delta_{ij}] \tag{4}$$

Where, $|V_i|$ and $|V_j|$ are voltage magnitudes of buses $i$ and $j$, and $\delta_{ij}$ is the phase difference between the two buses. Also, $\Delta G_{ij}$ and $\Delta B_{ij}$ can be expressed as below:

$$\Delta G_{ij} = \frac{x_{TCSC} r_{ij} (x_{TCSC} - 2x_{ij})}{(r_{ij}^2 + x_{ij}^2)[r_{ij}^2 + (x_{ij} - x_{TCSC})^2]} \tag{5}$$

$$\Delta B_{ij} = -\frac{x_{TCSC} r_{ij} (r_{ij}^2 - x_{ij}^2 + x_{TCSC} x_{ij})}{(r_{ij}^2 + x_{ij}^2)[r_{ij}^2 + (x_{ij} - x_{TCSC})^2]} \tag{6}$$

## 2.2 STATCOM

A FACTS device of the shunt kind called STATCOM may regulate bus voltage by absorbing or injecting reactive power [27]. By attaching to a suitable DC source, such as a battery or a fuel cell, it may also provide AC active power. Reactive and active power exchanges in the STATCOM structure will be regulated by the output voltage magnitude and phase regulation of the voltage source converter (VSC), respectively. It is possible to write the injection of STATCOM's active and reactive powers at the linked bus as follows:

$$P_i^S = \frac{1}{r_S^2 + x_S^2}\left[ r_S |V_i||V_S|\cos\delta_{iS} - x_S |V_i||V_S|\sin\delta_{iS} - r_S |V_i|^2 \right] \tag{7}$$

$$Q_i^S = \frac{1}{r_S^2 + x_S^2}\left[ r_S |V_i||V_S|\sin\delta_{iS} + x_S |V_i||V_S|\cos\delta_{iS} - x_S |V_i|^2 \right] \tag{8}$$

Where, $|V_S|, |V_i|, \delta_S, \delta_i$, are the voltage magnitudes and phase angles of shunt branch at bus $i$, respectively. Moreover, the operating limits of output voltage for STATCOM can be specified as:

$$|V_{S,\min}| \leq |V_S| \leq |V_{S,\max}| \tag{9}$$

Which, $\delta_S$ can take any angle in the range of $0$ and $2\pi$ [27].

# 3. DR Program Implementation

## 3.1 Concepts

Modern intelligent transmission networks have established a suitable framework for the positive involvement of the demand-side in the best possible performance of power systems. Demand side management (DSM) programs actually give customers the chance to reduce or move some of their time-shiftable loads during peak hours in response to financial incentives or changes in electricity prices, which could significantly increase flexibility in the safe and ideal operation of power systems. The incentive-based and price-based categories of the DR programs, which may be seen of as an authoritative CMP tool, were separated apart [6]. Power costs fluctuate over time under price-based schemes, which causes users' electricity usage habits to change [28]. In addition, the power provider offers incentives to consumers in the incentive-based DR program to reduce their electricity use when the system's security is in danger [29]. It should be noted that this article uses FACTS devices and an efficient incentive-based DR scheme to reduce congestion.

DR programs may be seen as effective CMP strategies, however DRPs alone can accomplish this goal in the best possible way. The best buses should be chosen based on ATC in order to execute the DR program effectively. The ISO develops an economic load model to create the incentive-based DR program in order to reduce congestion for the least amount of money. To this aim, the desired amount of demand at each bus of the network is determined using the price-elasticity concept connected with incentive and punishment policies. The sensitivity of demand in the ith period to the price in the same period is known as self-elasticity [29]:

$$E(i,i) = \frac{\rho_0(i)}{d_0(i)} \frac{\partial d(i)}{\partial \rho(i)} \tag{10}$$

Because there is no substitute for electrical energy in the near term, it is important to remember that in real-world operating difficulties, the price elasticity of demand is almost zero. Additionally, the cross-elasticity, which shows how sensitive the load is to the price in the *j*th period, may be described as:

$$E(i,j) = \frac{\rho_0(j)}{d_0(i)} \frac{\partial d(i)}{\partial \rho(j)} \quad (11)$$

Where, $d$ is the demand, $\rho$ is the price of elasticity and $E(i,i)$ as well as $E(i,j)$ are self and cross-elasticities, respectively. The subscript 0, stands for the condition before DR program implementation, i.e. before the change in the initial price.

As the elasticities of the consumers might differ from one to the next, a sizable database will be needed, which is often not accessible. As a result, the lack of data and precision in this kind of modeling might result in programming that is not exact. By modeling the price elasticity with regard to total demand, which will have the same value for all consumers, the issue of accuracy is overcome and a more realistic and useful DR technique is introduced in this study. In this kind of modeling, self and cross elasticities may be calculated from [29]-[31]:

$$E^P(i,i) = \frac{d_0(i)}{\rho_0(i)} \frac{\partial \rho(i)}{\partial d(i)} \quad (12)$$

$$E^P(i,j) = \frac{d_0(i)}{\rho_0(j)} \frac{\partial \rho(j)}{\partial d(i)} \quad (13)$$

Self-elasticity will be non-negative since if demand rises, the market price will rise or stay the same for ith duration. The cross elasticity will be negative if load shifting is available because rising demand in hour I will lead to falling demand in hour j and a falling or constant market price. The primary factor that affects how consumers behave in DR program models, such as time-of-use pricing, is elasticity. Elasticity and incentives play a major role in determining customer behavior in incentive-based DR systems like direct load management, and they also play a significant role in models like capacity markets [6]. As a result, while programming DR, these three important factors should be taken into account. The following is an one period DR program modeling taking the factors into account.

**3.2 Formulation**

In order to lower procurement costs and relieve line congestion, an incentive-based DR program called Interruptible/Curtailable (I/C) service has been presented in this study. Customers on I/C service tariffs get a discount or bill credit in return for agreeing to lower load when necessary under the planned I/C program, which is required. Customers will face consequences if they don't cut down. The following list of demand changes brought forth by the DR program:

International Transactions on Electrical Energy Systems$$\Delta d(i) = d(i) - d_0(i) \tag{14}$$

The customer benefit function ($S(i)$) is used to extract the final model of DR in this paper.

$$S(i) = B(d(i)) - d(i)\rho(i) \tag{15}$$

$$\frac{\partial S(i)}{\partial d(i)} = \frac{\partial B(d(i))}{\partial d(i)} - d(i)\frac{\partial \rho(i)}{\partial d(i)} - \rho(i) = 0 \tag{16}$$

$$\frac{\partial B(d(i))}{\partial d(i)} = d(i)\frac{\partial \rho(i)}{\partial d(i)} + \rho(i) \tag{17}$$

Substituting the initial values for customers demand in (17) will result in:

$$\frac{\partial B(d_0(i))}{\partial d(i)} = d_0(i)\frac{\partial \rho(i)}{\partial d(i)} + \rho_0(i) \tag{18}$$

On the other hand, the Taylor expansion of the revenue function can be written as:

$$B(d(i)) = B(d_0(i)) + \frac{\partial B(d_0(i))}{\partial d(i)}(d(i) - d_0(i)) \tag{19}$$

Deriving the above equation and using relations (12), (17) and (18):

$$\frac{\partial B(d(i))}{\partial d(i)} = \frac{\partial B(d_0(i))}{\partial d(i)} \tag{20}$$

$$\rho_0(i) - \rho(i) = (d(i) - d_0(i))\frac{\partial \rho(i)}{\partial d(i)} \tag{21}$$

$$d(i) = d_0(i)\left(1 + \frac{\rho_0(i) - \rho(i)}{\rho_0(i)E(i,i)}\right) \tag{22}$$

If the customer participates in an incentive-based DR program such as I/C program, $A(i)$ \$ will be paid as an incentive for each kWh of load reduction. So, the total incentive for participating in the I/C program can be obtained as:

$$P(\Delta d(i)) = A(i)[d_0(i) - d(i)] \tag{23}$$

On the other hand, if the participating customers in DR programs such as capacity market program do not comply with the obligations to reduce their load to the minimum acceptable level, they must pay the penalty which has been determined in the contract. If the penalty and requested load reduction level for *i*th period are denoted by $pen(i)$ and $IC(i)$ respectively, then the total penalty can be found as:

$$PEN(\Delta d(i)) = pen(i)\left[IC(i) - (d_0(i) - d(i))\right] \tag{24}$$

Rewriting the customer benefit function by considering the relations (23) and (24), single general period I/C program model will be obtained as:

$$d(i) = d_0(i)\left(1 + \frac{\rho_0(i) - \rho(i) - A(i) - pen(i)}{\rho_0(i)E(i,i)}\right) \tag{25}$$

Multi-period I/C program modelling, according to the proposed elasticity model can be found by extending the relation (25) for 24 hours of a day as below:

$$d(i) = d_0(i)\left(1 + \sum_{\substack{j=1 \\ j \neq i}}^{24} \frac{\rho_0(j) - \rho(j) - A(j) - pen(j)}{\rho_0(j)E(i,j)}\right) \tag{26}$$

Combining the relations (25) and (24) will give the final model of the I/C program, which is expressed in (27). This equation shows the optimal amount of participation for the enrolled customers in I/C programs. It is to be noted that the penalties and incentives will be set to zero for time-based programs.

$$d(i) = d_0(i)\left(1 + \frac{\rho_0(i) - \rho(i) - A(i) - pen(i)}{\rho_0(i)E(i,i)} + \sum_{\substack{j=1 \\ j \neq i}}^{24} \frac{\rho_0(j) - \rho(j) - A(j) - pen(j)}{\rho_0(j)E(i,j)}\right) \tag{27}$$

## 4. Problem formulation

This paper's primary goal is to offer an incremental welfare consensus method for wind integrated transmission networks based on the cooperative FACTS and DR program. A new two-stage market clearing method is created for this goal, with the first stage seeing GENCOs bid to the market with the intention of maximizing their personal profit and the second seeing ISO clear the market with the intention of maximizing the social welfare of the system. In the second stage, FACTS devices and the DR program are then handled effectively to relieve network congestion and lower GENCO's re-dispatch expenses. It is important to note that [32] is similar to the market clearing technique provided in this paper.

To determine the market clearing price in the day ahead pool based market, a two-stage optimization issue that takes into account numerous technological restrictions connected to generating units and transmission networks must be addressed. It is considered that every player is a taker of prices (i.e., there is no market power at the system). As a result, in order to participate in the market, all companies provide the ISO with their marginal costs. Equation (28), which has four parts, represents the day-ahead market clearing with the goal of maximizing the social welfare of the system. The first three parts express the offering of demands, the participation of GENCOs and wind power producers (WPP) in the market, and the fourth part shows the operation costs of GENCOs at each hour. This goal is constrained by four factors, including the power mismatch limitation (29) and the banned operating zones of market participants (30). (32).

$$Maximize: \sum_i (\lambda_i^D P_i^D) - \sum_j (\lambda_j^G P_j^G) - \sum_k (\lambda_k^W P_k^W) - \sum_j C_j(P_j^G) \tag{28}$$

$$\sum_i P_i^D - \sum_j P_j^G - \sum_k P_k^W = 0 \quad : \lambda^{DA} \tag{29}$$

$$0 \leq P_i^D \leq P_i^{\max} \qquad \forall i \in D \tag{30}$$

$$P_j^{\min} \leq P_j^G \leq P_j^{\max} \qquad \forall j \in G \tag{31}$$

$$P_k^{\min} \leq P_k^W \leq P_k^{\max} \qquad \forall k \in W \tag{32}$$

Where, $\lambda_i^D, \lambda_j^G, \lambda_k^W$ display the bidding prices of $i$th DISCO, $j$th GENCO and $k$th WPP that are announced to the ISO, respectively. Moreover, $P_i^D, P_j^G, P_k^W$ show the scheduled powers of $i$th demand, $j$th GENCO and $k$th WPP as well as $C(P_{gi})$ represents the operating cost function of GENCOs. The answer from (28) establishes each GENCO's and WPP's output as well as each customer's consumption at the same market clearing price. Each bus has a dual variable of active power balance constraint (29), which determines its LMP. Due to network congestion, the market pricing at each bus will have a different value. Three variables—the marginal cost of generation at the same bus, transmission congestion, and loss price—affect the precise value of LMP at each bus. It should be observed that the objective function (28) performs its computations notwithstanding any limits imposed by the communication network.

The CMP is implemented in the second stage as shown below [33]-[34] in order to take these restrictions into account while solving the stated issue. Equation (33) demonstrates the goal function of the second stage, which aims to reduce the cost associated with up- and down-power changes made by the ISO to relieve network congestion. It is important to note that any modification to the market-clearing circumstances will result in an additional payment to the participants whose output the ISO has modified. The first part in this equation indicates the entire cost of GENCOs' power adjustment for changes in power from their original schedule, while the second term states the cost borne by DR resources for reducing/increasing consumption.

$$\begin{aligned} Minimize: &\sum_j (r_j^{G,up} \Delta P_j^{G,up} + r_j^{G,down} \Delta P_j^{G,down}) \\ &+ \sum_i (r_i^{D,up} \Delta P_i^{D,up} x_i + r_i^{D,down} \Delta P_i^{D,down} y_i) \end{aligned} \tag{33}$$

Where, $\Delta P_j^{G,up}, \Delta P_j^{G,down}$ are increment and decrement of active power of $j$th GENCO in order to alleviate the congestion of network and also $r_j^{G,up}, r_j^{G,down}$ depict the offering prices of $j$th GENCO for increasing/decreasing its power for CMP purpose. $x_i, y_i$ are the DR commitment binary variables to increase and decrease consumption at any time which $x_i + y_i = 1$. Furthermore, $\Delta P_i^{D,up}, \Delta P_i^{D,down}$ and $r_i^{D,up}, r_i^{D,down}$, denote the International Transactions on Electrical Energy Systemsconsumption alterations and offered prices of *i*th DISCO for increasing/ decreasing its demand for CMP, in which $r_i^{D,up}, r_i^{D,down}$ are ramping limitations of the demand. It should be stated that here the AC power flow is used because the accuracy of this approach is better than DC one [35].

Due to a variety of technological limitations brought on by components and the network, the second stage is reduced. The active and reactive balance at each bus of the system is confirmed by equations (34) and (35) Each bus's nodal active power is equal to the total of the original demand plus the decrease after the deployment of DR. Additionally, the power factor and active power of the load are used to determine the reactive power of each load node. Furthermore, during an operation horizon time, restrictions (36) to (41) enforce the participants' permissible boundaries. The ramp-up/down restrictions of GENCOs are shown in constraints (42) and (43), and the combined powers of DISCOs and GENCOs are shown in constraints (44) and (45). Note that $P_i^{D,DA}, P_j^{G,DA}$ in these equations are obtained from the first stage, therefore, in the stage they are fixed to their optimal value. The permissible bound of voltage is enforced by (46) as well as the active and reactive powers passing through the transmission lines are limited by constraints (47) to (48). In the end, equations (49) and (48) illustrate the active and reactive power flow equations, respectively.

$$P_n^G + P_n^W - P_n^D = p_n(V,\theta) \qquad \forall n \in \mathbb{N} \qquad (34)$$

$$Q_n^G + Q_n^W - Q_n^D = q_n(V,\theta) \qquad \forall n \in \mathbb{N} \qquad (35)$$

$$P_i^{\min} \leq P_i^D \leq P_i^{\max} \qquad \forall i \in D \qquad (36)$$

$$P_j^{\min} \leq P_j^G \leq P_j^{\max} \qquad \forall j \in G \qquad (37)$$

$$P_k^{\min} \leq P_k^W \leq P_k^{\max} \qquad \forall k \in W \qquad (38)$$

$$Q_i^D = P_i^D \tan(\phi_i^D) \qquad \forall i \in D \qquad (39)$$

$$\sqrt{(S_j^G)^2 - (P_j^{\min})^2} \leq Q_j^G \leq \sqrt{(S_j^G)^2 - (P_j^{\max})^2} \qquad \forall j \in G \qquad (40)$$

$$-\frac{P_k^W \sqrt{1-(\kappa_k^{\min})^2}}{\kappa_k^{\min}} \leq Q_k^W \leq \frac{P_k^W \sqrt{1-(\kappa_k^{\min})^2}}{\kappa_k^{\min}} \qquad \forall k \in W \qquad (41)$$

$$P_{j,t}^G - P_{j,t-1}^G \leq RU_j \qquad \forall j \in G \qquad (42)$$

$$P_{j,t-1}^G - P_{j,t}^G \leq RD_j \qquad \forall j \in G \qquad (43)$$

$$P_i^D = P_i^{D,DA} + \Delta P_i^{D,up} - \Delta P_i^{D,down} \qquad \forall i \in D \qquad (44)$$

$$P_j^G = P_j^{G,DA} + \Delta P_j^{G,up} - \Delta P_j^{G,down} \qquad \forall j \in G \qquad (45)$$

$$V_n^{\min} \leq V_n \leq V_n^{\max} \qquad \forall n \in \mathbb{N} \qquad (46)$$

$$P_b \leq P_b^{\max} \qquad \forall b \in B \qquad (47)$$

$$Q_b \leq Q_b^{\max} \qquad \forall b \in B \qquad (48)$$

$$P_b = G_b(\cos(\delta_n - \delta_m) - 1) + B_b \sin(\delta_n - \delta_m) \qquad \forall b \in B \qquad (49)$$

$$Q_b = G_b(\sin(\delta_n - \delta_m) - 1) - B_b \cos(\delta_n - \delta_m) \qquad \forall b \in B \qquad (50)$$

In these equations, $P_n^G, P_n^W, P_n^D$ illustrate total final active powers of GENCO, WPP, and DISCO at the $n$th bus, respectively, and similarly $Q_n^G, Q_n^W, Q_n^D$ are total final reactive powers generated by GENCO, WPP, and DISCO at the $n$th bus; $p_n(V,\theta), q_n(V,\theta)$ express the injection of active and reactive powers in the $n$th bus with respect to magnitude and angle of its voltage. $\kappa_k^{\min}$ is the minimum power factor of WT which can be adjusted by controlling the inverters (i.e., 0.9 lagging, 0.9 leading) and $RU_j, RD_j$ refer to ramp up/down of $j$th GENCO for increment/decrement of its production. Moreover, $V_n$ shows the voltage magnitude of $n$th bus of system, $\delta_n$ is the voltage angle of $n$th bus, $P_b, Q_b$ display the active and reactive powers crossed from $b$th line and $G_b, B_b$ are conductance and susceptance of $b$th line, respectively.

Additionally, two separate kinds of FACTS devices, TCSC and STATCOM, are used to lower re-dispatch costs and help the ISO relieve network congestion. By properly using these controllers, congestion reduction, LMP smoothing, and voltage security enhancement may all be achieved. These FACTS controllers' operational limitations may be represented as:

$$x_{TCSC,\min} \leq x_{TCSC} \leq x_{TCSC,\max} \qquad (51)$$

$$B_{STAT}^{\min} \leq B_{STAT} \leq B_{STAT}^{\max} \qquad (52)$$

As mentioned before, the TCSC can be modeled as a changeable reactance ($x_{TCSC}$) during the steady-state condition as (51). As the same way, an STATCOM is modeled as a variable susceptance ($B_{STAT}$) as shown (52) by voltage source converter with a specified range of voltage magnitude.

Fig. 1 illustrates the proposed two-stage algorithm for CMP of wind integrated transmission networks. This paper's major contribution is to lay forth a framework for coordination between FACTS devices and DR algorithms in order to provide congestion control at the lowest feasible cost. As is evident, the GENCOs and DISCOs offer their pricing and power to the market, and the ISO then clears the market based on the proposals it has received from participants. The ISO then assesses the technical limitations of the system and, in the event of line congestion, employs the FACTS devices to ease the situation as far as is practical. If the FACTS devices fail to relieve the congestion, the ISO implements the DR to lower consumption and so relieve the congestion.

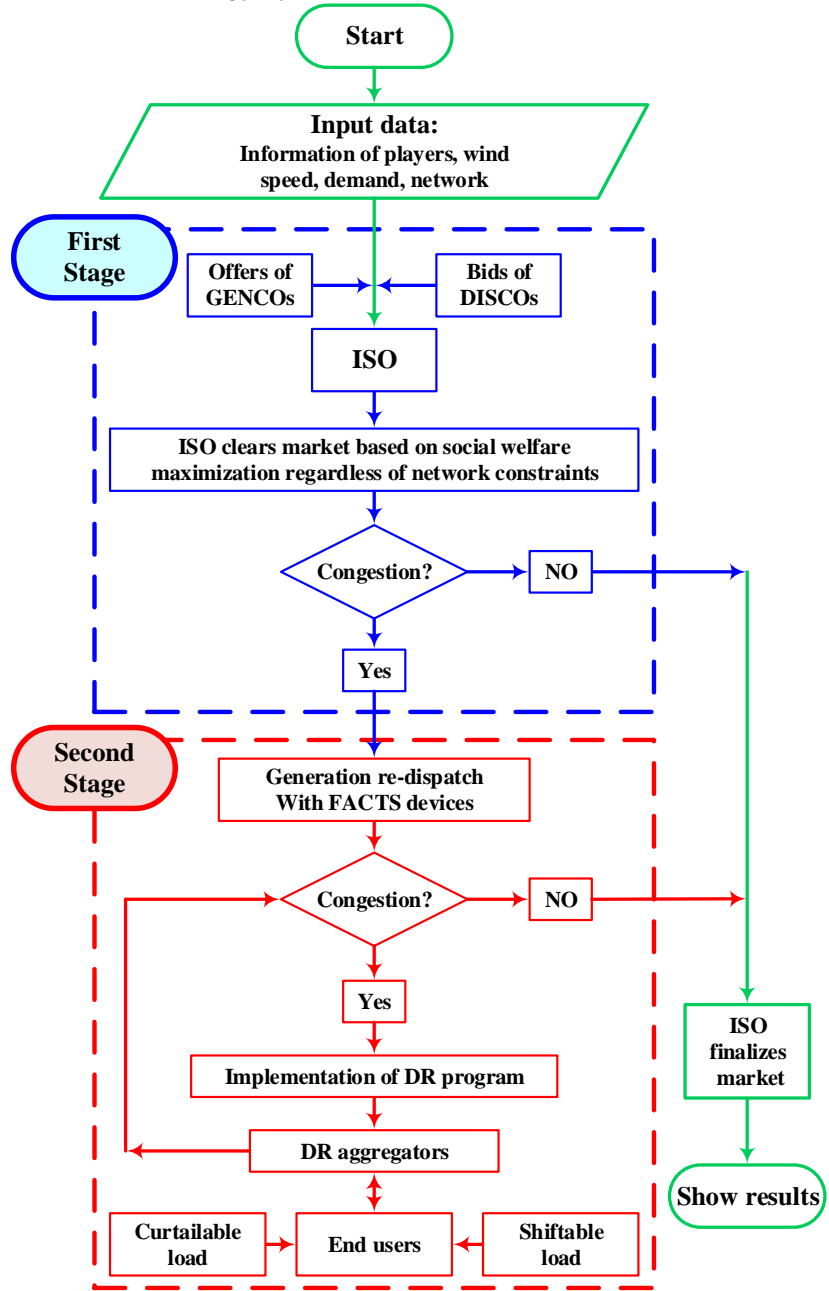

**Fig. 1.** Flowchart of proposed integrated method for CMP of transmission networks.

# 5. Simulation results

## 5.1 Data and case study

The operators must have powerful tools to offset their unpredictability and uncertainty given the fast integration of WPPs into the power networks. The promotion of energy storage technology has historically been used to solve this issue. To control the security and stability of the market, however, the ISO might profit from other choices inside the electrical markets, such as DSM schemes [36]. The key issue is that WPP intermittent periods

International Transactions on Electrical Energy Systems

occur during peak load hours; as a result, DSM technologies may be used to lower consumption during these times [37]. In order to relieve the congestion of wind integrated transmission networks in the liberalized power market, this research suggests a joint FACTS and DR scheme. This is done by using a real-world case study based on the modified IEEE 24-bus reliability test system (RTS), as illustrated in Fig. 2 [38], to confirm the effectiveness and utility of the suggested strategy.

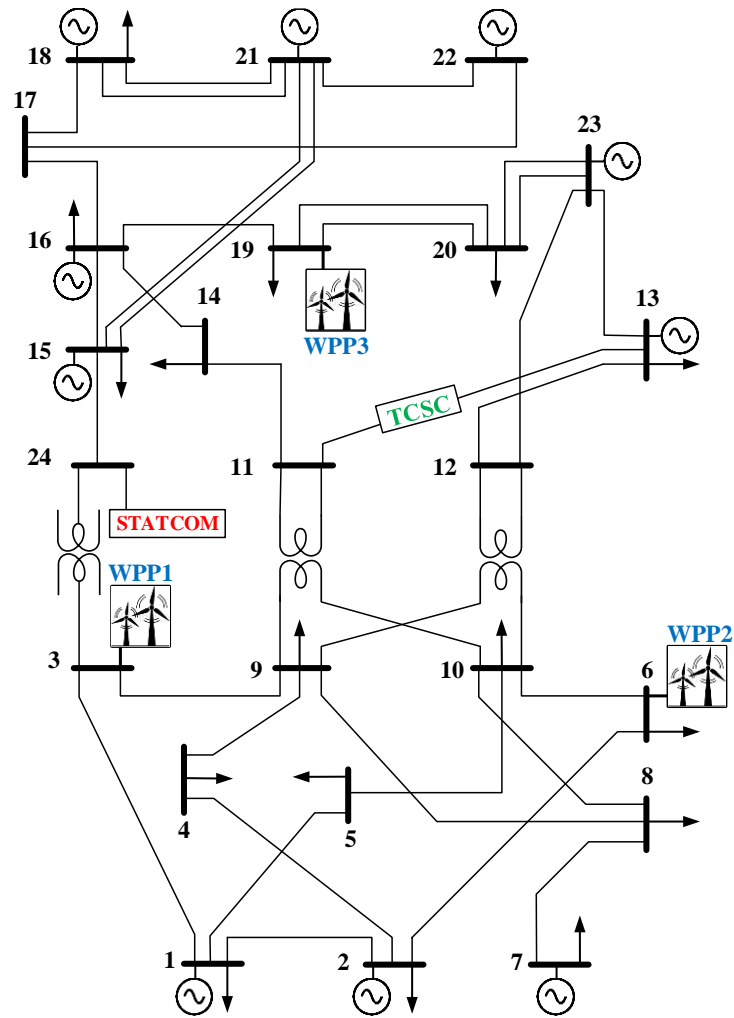

**Fig. 2.** IEEE one-area 24 bus RTS with the integration of WPP and FACTS devices.

The network's generators are organized simply by kind and node. Table 1 provides the self and cross elasticity of demand, which is needed to carry out the suggested DR program. Additionally, Table 2 shows the candidate bus and amount of DR resources to take part in the CMP. The suggested model makes the assumption that only 10 buses may participate in the DR as curtailable loads, with a 30% demand participation ratio. It should be noted that the penalty and incentive rates are set at 142 and 85 dollars per megawatt hour, respectively. Table 3 displays the information on the reactance/susceptance limits of the FACTS devices that have been used. In reality, these FACTS devices may change the pattern of power flow at the network in such a manner as to overcome the congestion in overloaded lines by modifying the reactance and susceptance of lines. The quadratic model is used in this study to simulate the operational expenses of GENCOs, whose cost functions are shown in Table 4. The availability of GENCOs' offerings is based on their cost functions and available capacities, which may be computed as follows.

$$C_j(P_j^G) = \alpha_j + \beta_j P_j + \gamma_j P_j^2 \qquad \forall j \in G \qquad (65)$$

**Table 1.** Self and cross elasticity of demand

|  | Peak | Off-peak | Low |
|---|---|---|---|
| Peak | -0.12 | 0.018 | 0.014 |
| Off-peak | 0.018 | -0.10 | 0.012 |
| Low | 0.014 | 0.012 | -0.10 |

**Table 2.** Candidate bus and size of DR resources

| DR number | Site (Bus Number) | Size (MW) |
|---|---|---|
| 1 | 3 | 64.79 |
| 2 | 6 | 54.82 |
| 3 | 8 | 69.04 |
| 4 | 9 | 73.22 |
| 5 | 10 | 76.48 |
| 6 | 13 | 95.73 |
| 7 | 15 | 127.46 |
| 8 | 16 | 45.79 |
| 9 | 18 | 148.25 |
| 10 | 20 | 65.71 |

**Table 3.** Data of FACTS devices

| Type of FACTS | TCSC | STATCOM |
|---|---|---|
| Operating limit | $-0.25 \leq X_{TCSC} \leq 0.25$ | $-0.33 \leq B_{STAT} \leq 0.33$ |
| Location | Line 18 | Bus 24 |

Three private WPPs are taken into account in the proposed model and are situated at buses 3, 6, and 19. Table 5 lists some of their features. Additionally, the PEM technique is used to simulate the uncertainty of wind production. It should be underlined that WPPs should be used to their full potential in terms of economic concerns. In other words, the use of FACTS devices should be used to limit the wind power leakage of system, which arises for technical reasons such as inadequate transmission capacity. WPPs have not participated in the CMP because of this.

Table 4. Cost coefficients of generators

| GENCO number | α | β | γ |
|---|---|---|---|
| 1 | 212.3076 | 16.0811 | 0.014142 |
| 2 | 212.3076 | 16.0811 | 0.014142 |
| 3 | 781.5212 | 43.6615 | 0.052672 |
| 4 | 832.7575 | 48.5804 | 0.007170 |
| 5 | 86.38524 | 56.5646 | 0.328412 |
| 6 | 382.2391 | 12.3883 | 0.008342 |
| 7 | 395.3749 | 4.42317 | 0.000213 |
| 8 | 395.3749 | 4.42317 | 0.000213 |
| 9 | 382.2391 | 12.3883 | 0.008342 |
| 10 | 665.1094 | 11.8495 | 0.004895 |

Table 5. Characteristics of WTs

| Number | $P_{min}$ (MW) | $P_{max}$ (MW) | $V_r$ (m/s) | $V_{cut-in}$ (m/s) | $V_{cut-out}$ (m/s) |
|---|---|---|---|---|---|
| WPP1 | 50 | 110 | 12 | 8 | 20 |
| WPP2 | 60 | 125 | 16 | 10 | 24 |
| WPP3 | 60 | 125 | 16 | 10 | 24 |

## 5.2 Numerical results

The proposed problem is a mixed-integer non-linear program (MINLP) which is solved by DICOPT solver under GAMS software [39] and implemented on a laptop with Intel Core i7 3.50 GHz CPU and 32 GB RAM. The average implementation time of the proposed algorithm is 93.84 seconds. It should be mentioned that both the absolute gap and relative gap are set to be zero (i.e., optca=0, optcr=0).

Table 6. Results of market clearing process (first stage)

| Gen. no | Generation (MW) | Demand. no | Consumption (MW) |
|---|---|---|---|
| 1 | 178.64 | 1 | 197.64 |
| 2 | 178.64 | 2 | 181.37 |
| 3 | 205.52 | 3 | 236.59 |
| 4 | 215.88 | 4 | 249.33 |
| 5 | 149.30 | 5 | 187.05 |
| 6 | 137.14 | 6 | 341.62 |
| 7 | 374.96 | 7 | 228.70 |
| 8 | 366.48 | 8 | 294.44 |
| 9 | 285.07 | 9 | 325.11 |
| 10 | 541.95 | 10 | 371.58 |

**Table 7.** Results of generation re-dispatch in different cases

| Gen. no | Without FACTS$DR | Only with DR | Only with FACTS | With FACTS$DR |
|---|---|---|---|---|
| 1 | 16.88 | 10.54 | 13.25 | 7.23 |
| 2 | 15.42 | 8.76 | 14.36 | 6.29 |
| 3 | 3.64 | 3.52 | 3.65 | 1.05 |
| 4 | 9.81 | 5.89 | 8.94 | 2.67 |
| 5 | -0.74 | 0.02 | 2.14 | -0.05 |
| 6 | -5.23 | -2.28 | -3.48 | -0.98 |
| 7 | -17.65 | -15.47 | -18.21 | -6.80 |
| 8 | 8.62 | 8.01 | 6.74 | 2.84 |
| 9 | -11.25 | -9.34 | -11.04 | -6.66 |
| 10 | 12.47 | 10.22 | 12.39 | 5.01 |

**Table 8.** The setting of FACTS devices

| FACTS | Reference setting (p.u) |
|---|---|
| TCSC | -0.0923 |
| STATCOM | 0.0746 |

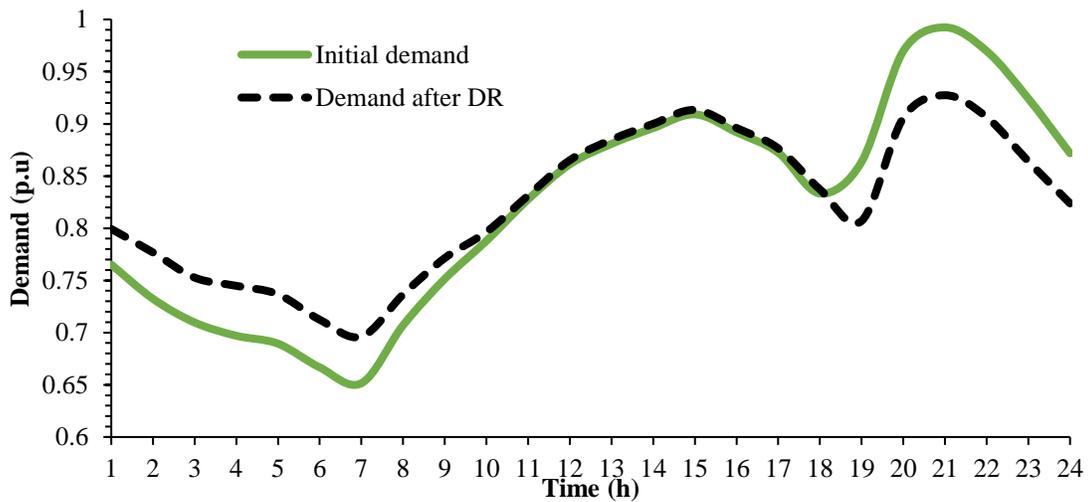

**Fig. 3.** Daily load profile of the system with and without DR applying.

Table 6 displays the outcomes of the first step, which included market clearing despite congestion. Following market clearing, the ISO should assess network congestion before implementing DR and rescheduling generators to help the system become less congested. To demonstrate the effectiveness of the suggested integrated strategy, the results of the second stage taking various scenarios into consideration are shown in Tables 7 and 8. Additionally, Fig. 6 shows the load profile of the system with and without the use of the DR software. Table 7 reveals that the overall number of generators re-dispatching in the simultaneous FACTS and DR implementation is much less than in the other three situations (positive value is to increase output power and negative values is to reduce it). This problem demonstrates that the joint operation of FACTS devices and DR program will provide better outcomes at a cheaper cost than each program operating alone. In other words, these two tools work best together to provide outcomes that are more palatable. The DR program moves a portion of the system's peak load to off-peak times (as shown in Fig. 3), and FACTS devices change the characteristics of the lines to change the network's power flow pattern in order to effectively reduce line congestion. Various scenarios of the electricity flow of lines during peak hours are shown in Fig. 4. As can be observed, by implementing the FACTS and DR programs, the volume of power transactions in the overloaded lines is significantly decreased. This is because customers' involvement in the DR program to cut down on their usage has distorted the flow of electricity to low load lines. In reality, FACTS devices will cause power to be diverted to low-load lines in order to balance the power transactions of the lines via the optimum configuration of line reactance. Another conclusion drawn from this graph is that the suggested integrated strategy, as opposed to using each component separately, has a better effect on the system's CMP.

The development of various nodal pricing throughout the network is one of the key effects of congestion in transmission networks. LMP derives from market clearing for all players, which enables energy exchange at each network bus to be completed based on the nodal energy price in the same bus. For the purpose of managing both FACTS device and DR resources, the proposed CMP algorithm uses the LMP as an input. Fig. 5 shows the LMP at each bus in the network under various scenarios. It is clear that the suggested integrated solution, which is based on FACTS and DR participating as effectively as possible, makes the energy cost more uniform than in other scenarios at various buses. In other words, the suggested algorithm makes the system's social welfare higher than it would be if these resources were used independently.

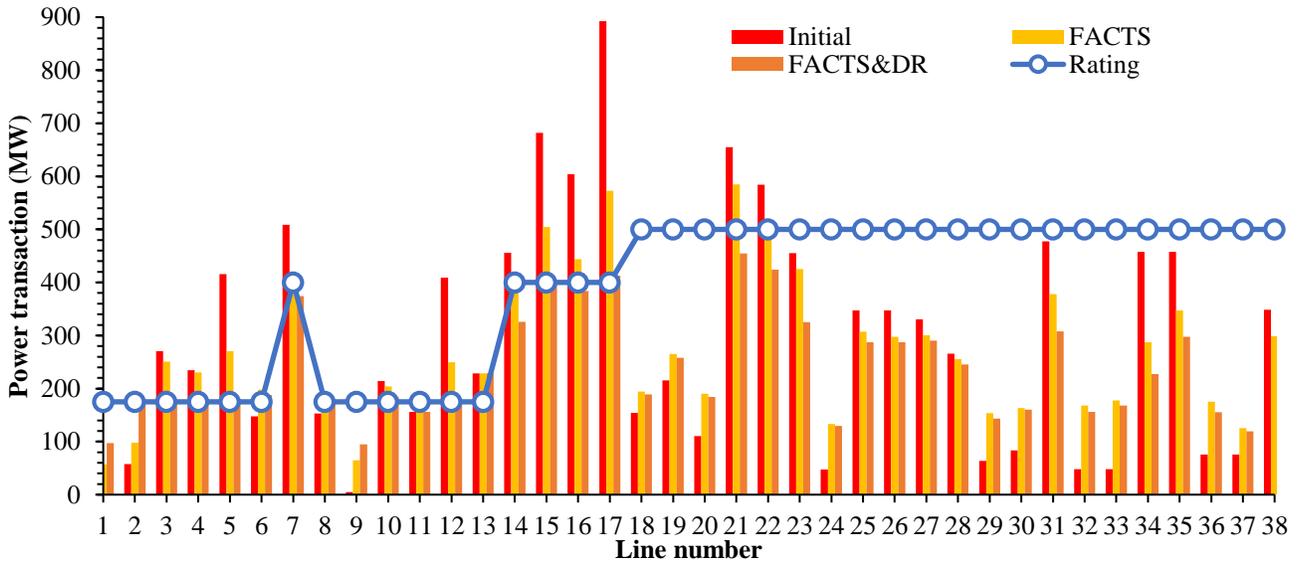

**Fig. 4.** Transaction powers through the lines in different cases.

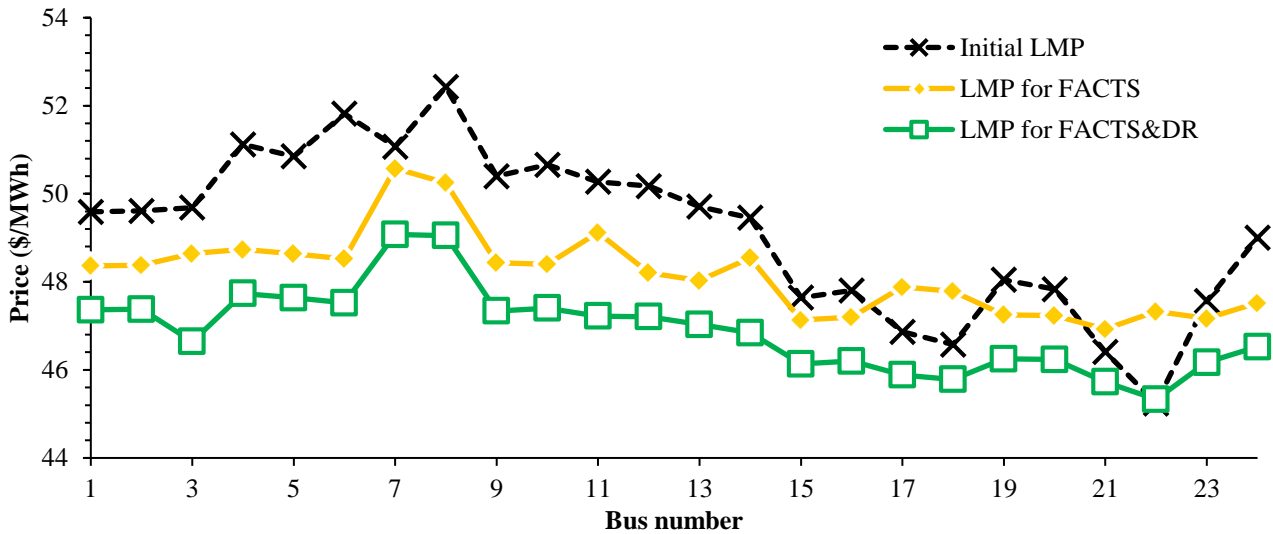

**Fig. 5.** Value of LMP at different buses of network with and without applying of the proposed method.

Fig. 6 shows the daily delivery of GENCOs, which was collected from the first stage. To limit the system's income sufficiency, GENCOs are sent in a merit-order manner. To get the most profit at the lowest cost, WPPs are operated at their maximum position in the suggested algorithm. Additionally, Figs. 7 and 8 demonstrate that the suggested method not only improves system technical parameters like power losses and voltage profile, but also lessens network congestion. In addition, the power losses and voltage profiles of the system will be addressed by enhancing the power flow pattern in the transmission network and directing surplus power in high-loaded lines towards low-loaded lines.

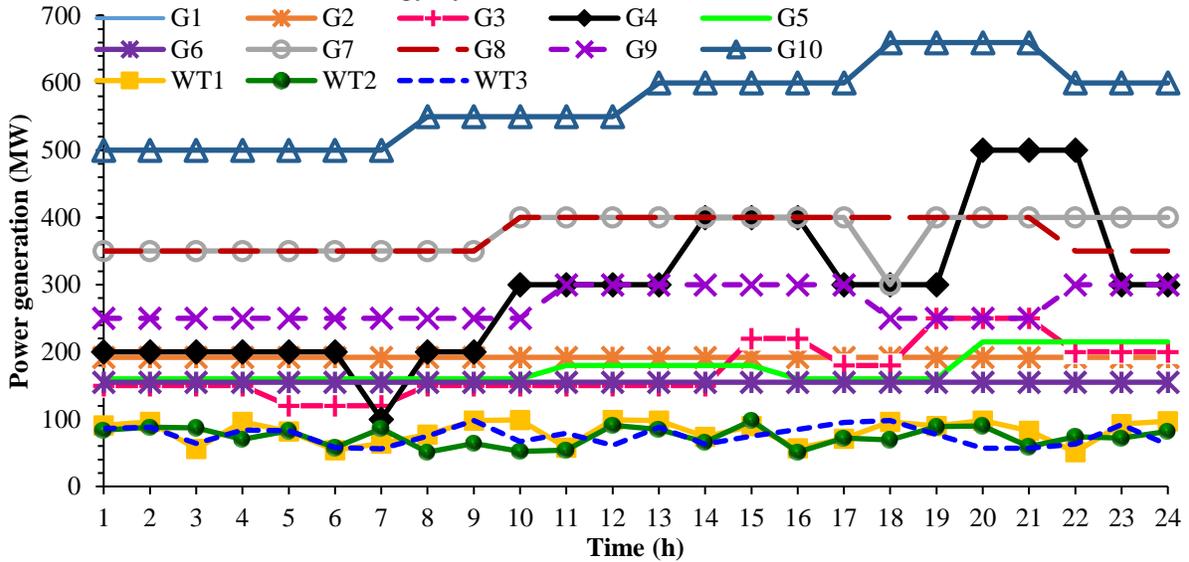

**Fig. 6.** Power generations of various GENCOs and WPPs over entire operation horizon.

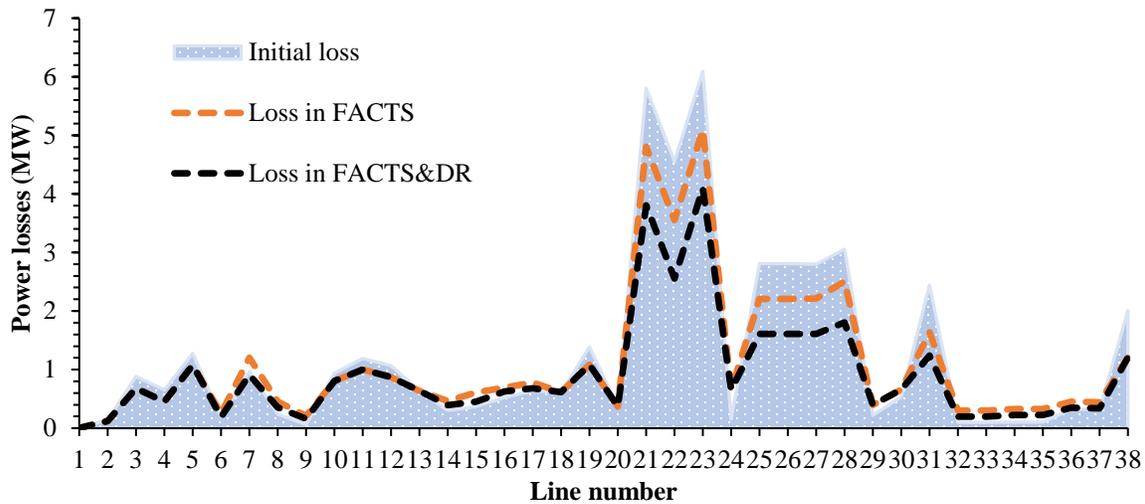

**Fig. 7.** Power losses of lines in different cases.

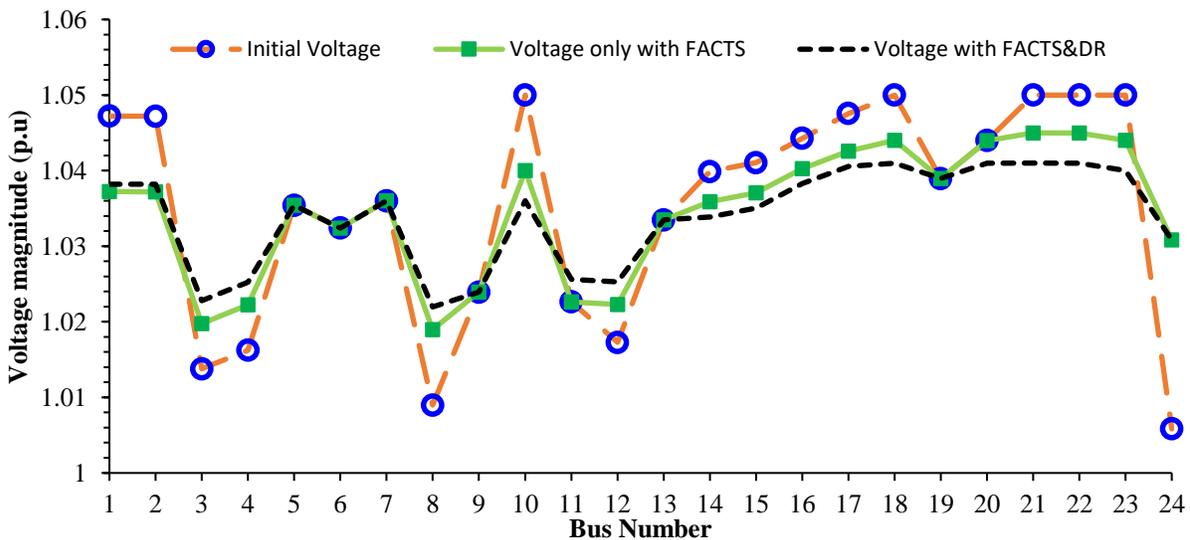

**Fig. 8.** Voltage profile of system in different cases.

# 6. Conclusion

In this research, an unique algorithm for the CMP of wind-integrated transmission networks was presented. This algorithm augments generation re-dispatch and FACTS devices with DR implementation that is dependent on market mechanisms. By doing this, an unique two-stage optimization was created in order to minimize costs while maximizing the social welfare of the whole system and reducing network congestion. The LMP is used as an input by the proposed CMP algorithm to reward DR resources and FACTS devices for reducing network congestion and provide voltage support. In addition, the suggested method finds the best location and quantity of DR resources to control system congestion by varying electricity pricing tariffs over time in an effort to persuade consumers to alter their energy consumption habits.

The numerical results from simulations have demonstrated that the distributions of all state variables and line-flow quantities can be accurately and effectively evaluated with the proposed PEM approach through straightforward numerical computations, provided the uncertain parameters can be measured or estimated. However, the suggested approach might employ more estimatepoints for improved outcomes. The suggested method is useful for large-scale case studies since it is precise and doesn't need a lot of processing effort. Additionally, it has been shown that combining the FACTS and DR programs may lower congestion costs more effectively than doing so with only one of them. Additionally, the suggested solution may be shown to address the technical problems with the system, such as power losses and voltage stability, in addition to reducing system congestion.